\documentclass[a4paper,12pt]{article}
\usepackage{amsmath}
\usepackage{latexsym,amssymb,amsbsy,natbib,rotating,setspace,amsxtra}
\usepackage[T1]{fontenc}
\usepackage[latin1]{inputenc}
\usepackage[english]{babel}
\usepackage{graphicx,pict2e}
\usepackage{algorithm}
\floatname{algorithm}{Algorithm}

\accentedsymbol{\betahathat}{\Hat{\Hat \beta}}

\newtheorem{theo}{Theorem}
\newtheorem{theoproof}{Proof.}

\newtheorem{prop}{Proposition}
\newtheorem{propproof}{Proof.}

\newtheorem{lemma}{Lemma}
\newtheorem{lemproof}{Proof.}

\title{Martingale central-limit theorems for pivotal sampling
~\\ ~ \\
\large October 29th, 2015}

\date{}

\author{Guillaume Chauvet
\thanks{Laboratoire de Statistique d'Enqu\^ete, CREST/ENSAI, Campus de Ker Lann, 35170 Bruz, France}}

\begin{document}

\maketitle

\begin{abstract}
\noindent Ordered pivotal sampling is one of the simplest algorithm to perform without-replacement unequal probability sampling. It has found uses in the context of longitudinal surveys and spatial sampling, and enables in particular a good spatial balance of the selected units. In this work, we follow the approach proposed by Ohlsson~(1986), and apply a martingale central-limit theorem to prove the asymptotic normality of the Horvitz-Thompson estimator under a design-based approach, and under a model-assisted approach. In particular, our model assumptions allow for correlations between values, which is of particular interest for applications in spatial sampling.
\end{abstract}


\normalsize

\section{Introduction} \label{sec:intro}

\noindent Many different procedures exist in the finite population context for without-replacement unequal probability sampling. The ordered pivotal sampling algorithm proposed by Deville and Till\'e (1998) and rediscovered by Srinivasan (2001) is one of the simplest procedures. Some important properties of this algorithm have been derived recently, see for example Br\"anden and Jonasson (2011), Chauvet (2012) and Jonasson (2012), and the references therein. Pivotal sampling has found uses in the context of longitudinal surveys (Nedyalkova et al., 2009) and spatial sampling (Grafstrom et al., 2012). In case of a multidimensional population, these last authors point out that ordered pivotal sampling enables to select samples well spread over space, which may result in more efficient estimators for totals if the variable of interest presents a positive spatial correlation. \\

\noindent So as to supply estimates with confidence intervals with appropriate coverage, it is of interest to state a central-limit theorem for the Horvitz-Thompson estimator. In case of pivotal sampling, this result would follow from Proposition 1 in Br\"anden and Jonasson (2011), but their assumptions require in particular that the variance of the Horvitz-Thompson estimator under pivotal sampling is asymptotically equivalent to that under Poisson sampling, which is not straightforward to prove. In the current paper, we first prove the asymptotic normality of the Horvitz-Thompson estimator for pivotal sampling under the randomization associated to the sampling design. We follow the approach in Ohlsson~(1986) for this purpose, and apply a martingale central-limit theorem. Since some second-order inclusion probabilities are usually equal to zero, customary design-based variance estimators are usually biased. We therefore appeal to some super-population model for the variable of interest. The asymptotic normality of the Horvitz-Thompson estimator is proved under the randomization associated to both the sampling design and the super-population model, under some mild conditions. In particular, our model assumptions allow for correlations between values, which is of particular interest for applications in spatial sampling. \\

\noindent The paper is organized as follows. In Section \ref{sec:nota}, the notation is defined. A recursive algorithm for ordered pivotal sampling is presented in Section \ref{sec:ops}, and some useful properties are derived. They are used in Section \ref{sec:tcl} to obtain a central-limit theorem for the Horvitz-Thompson estimator under the randomization associated to the sampling design. In Section \ref{sec:model}, we appeal to some super-population model for the variable of interest, and the asymptotic normality of the Horvitz-Thompson estimator is proved under the randomization associated to both the sampling design and the super-population model. We make some final remarks in Section \ref{sec:conc}. Some additional lemma are given in Appendix.

\section{Notation} \label{sec:nota}

\noindent We consider a finite population $U$ of size $N$. In order to study the asymptotic properties of the sampling designs and estimators that we treat below, we consider the asymptotic framework of Isaki and Fuller (1982). We assume that the population belongs to a nested sequence $\{U_{\nu}\}$ of finite populations with increasing sizes $N_{\nu}$, and all limiting processes will be taken as $\nu \to \infty$. Though all quantities under consideration depend on $\nu$, this subscript is omitted in what follows for simplicity of notation. \\

\noindent Denote $\pi_U=\left(\pi_{1},\ldots,\pi_{N}\right)^{\top}$ a vector of probabilities, with $0 < \pi_{k} \leq 1$ for any unit $k$ in $U$ and $n=\sum_{k \in U} \pi_k$ the sample size. Also, denote
    \begin{eqnarray} \label{max:pik}
      \pi_M & = & \max_{k \in U} \pi_k.
    \end{eqnarray}
We are interested in estimating the total $t_{y}=\sum_{k \in U} y_{k}$ for some variable of interest taking the value $y_{k}$ for unit $k \in U$. We note $y_U=(y_1,\ldots,y_N)^{\top}$ the vector of the population values. In Sections \ref{sec:ops} and \ref{sec:tcl}, $y_U$ is seen as a vector of deterministic quantities and a central-limit theorem is obtained under the sole randomization associated to the sampling design. In Section \ref{sec:model}, we appeal to some super-population model. The vector $y_U$ is then seen as a random vector, and inference is made with respect to the randomization associated to both the sampling design and the super-population model. \\

\noindent Denote $S$ a random sample, selected by means of a sampling design $p(\cdot)$ with inclusion probabilities $\pi_U$. Then the Horvitz-Thompson (HT) estimator
    \begin{eqnarray}
      \hat{t}_{y\pi} & = & \sum_{k \in S} \check{y}_k
    \end{eqnarray}
is design-unbiased for $t_y$, with $\check{y}_k=y_k/\pi_k$. We note $E(\cdot)$ and $V(\cdot)$ for the expectation and the variance of some estimator. Also, we note $E_{\{X\}}(\cdot)$ and $V_{\{X\}}(\cdot)$ for the expectation and variance conditionally on some random variable $X$.\\

\noindent We define $V_{k}=\sum_{l=1}^k \pi_{l}$ for any unit $k \in U$, and $V_0=0$. A unit $k \in U$ is said to be {\em cross-border} if $V_{k-1} < i$ and $V_{k} \ge i$ for some positive integer $i$. The cross-border units are denoted as $k_i,~i=0,\ldots,n$. We note
    \begin{eqnarray*}
      a_i=i-V_{k_i-1} & \textrm{ and } & b_i=V_{k_i}-i
    \end{eqnarray*}
for $i=1,\ldots,n-1$. For $k_0$ and $k_n$, we take by convention $a_0=b_0=0$ and $a_n=b_n=0$. The units $k_0$ and $k_n$ are in fact phantom units with zero associated probabilities. \\

\begin{figure}[htb!]
    \setlength{\unitlength}{0.03cm}
    \begin{picture}(600,200)
    \thicklines
    \put(30,130){\line(1,0){400}} \put(20,130){\line(1,0){5}} \put(10,130){\line(1,0){5}} \put(0,130){\line(1,0){5}} \put(435,130){\line(1,0){5}} \put(445,130){\line(1,0){5}} \put(455,130){\line(1,0){5}}
    \put(130,120){\line(0,1){20}} \put(130,145){\makebox(0,0)[b]{$i-1$}} \put(330,120){\line(0,1){20}} \put(330,145){\makebox(0,0)[b]{$i$}}
    \put(90,125){\line(0,1){10}}   \put(90,110){\makebox(0,0)[b]{\tiny{$V_{k_{i-1}-1}$}}} \put(160,125){\line(0,1){10}}   \put(160,110){\makebox(0,0)[b]{\tiny{$V_{k_{i-1}}$}}}
    \put(210,125){\line(0,1){10}}   \put(210,110){\makebox(0,0)[b]{\tiny{$V_{k_{i-1}+1}$}}} \put(240,115){\makebox(0,0)[b]{\tiny{$\ldots$}}}
    \put(270,125){\line(0,1){10}}   \put(270,110){\makebox(0,0)[b]{\tiny{$V_{k_{i}-2}$}}} \put(300,125){\line(0,1){10}}   \put(300,110){\makebox(0,0)[b]{\tiny{$V_{k_{i}-1}$}}}
    \put(350,125){\line(0,1){10}}   \put(350,110){\makebox(0,0)[b]{\tiny{$V_{k_{i}}$}}} \put(390,125){\line(0,1){10}}   \put(390,110){\makebox(0,0)[b]{\tiny{$V_{k_{i}+1}$}}}
    \put(90,85){\vector(1,0){38}}  \put(128,85){\vector(-1,0){38}} \put(110,75){\makebox(0,0)[b]{\tiny{$a_{i-1}$}}} \put(130,85){\vector(1,0){28}} \put(158,85){\vector(-1,0){28}} \put(145,75){\makebox(0,0)[b]{\tiny{$b_{i-1}$}}}
    \put(300,85){\vector(1,0){28}} \put(328,85){\vector(-1,0){28}} \put(315,75){\makebox(0,0)[b]{\tiny{$a_{i}$}}} \put(330,85){\vector(1,0){18}} \put(348,85){\vector(-1,0){18}} \put(340,75){\makebox(0,0)[b]{\tiny{$b_{i}$}}}
    \put(90,170){\vector(1,0){68}}  \put(158,170){\vector(-1,0){68}} \put(125,175){\makebox(0,0)[b]{\tiny{$\pi_{k_{i-1}}$}}} \put(160,170){\vector(1,0){48}}  \put(208,170){\vector(-1,0){48}} \put(185,175){\makebox(0,0)[b]{\tiny{$\pi_{k_{i-1}+1}$}}} \put(235,180){\makebox(0,0)[b]{\tiny{$\ldots$}}}
    \put(270,170){\vector(1,0){28}}  \put(298,170){\vector(-1,0){28}} \put(285,175){\makebox(0,0)[b]{\tiny{$\pi_{k_{i}-1}$}}} \put(300,170){\vector(1,0){48}}  \put(348,170){\vector(-1,0){48}} \put(325,175){\makebox(0,0)[b]{\tiny{$\pi_{k_{i}}$}}} \put(350,170){\vector(1,0){38}}  \put(388,170){\vector(-1,0){38}} \put(375,175){\makebox(0,0)[b]{\tiny{$\pi_{k_{i}+1}$}}}
    \put(30,65){\line(1,0){100}} \put(130,65){\line(0,1){10}} \put(65,50){\makebox(0,0)[b]{$U_{i-1}$}}
    \put(20,65){\line(1,0){5}} \put(10,65){\line(1,0){5}} \put(0,65){\line(1,0){5}}
    \put(130,30){\line(1,0){200}} \put(130,30){\line(0,1){10}} \put(330,30){\line(0,1){10}} \put(220,15){\makebox(0,0)[b]{$U_{i}$}}
    \put(330,65){\line(1,0){100}} \put(330,65){\line(0,1){10}} \put(405,50){\makebox(0,0)[b]{$U_{i+1}$}}
    \put(435,65){\line(1,0){5}} \put(445,65){\line(1,0){5}} \put(455,65){\line(1,0){5}}
    \end{picture}
    \vspace*{-1cm} \caption{Probabilities and cross-border units in microstratum $U_i$, for population $U$} \label{nota}
    \end{figure}
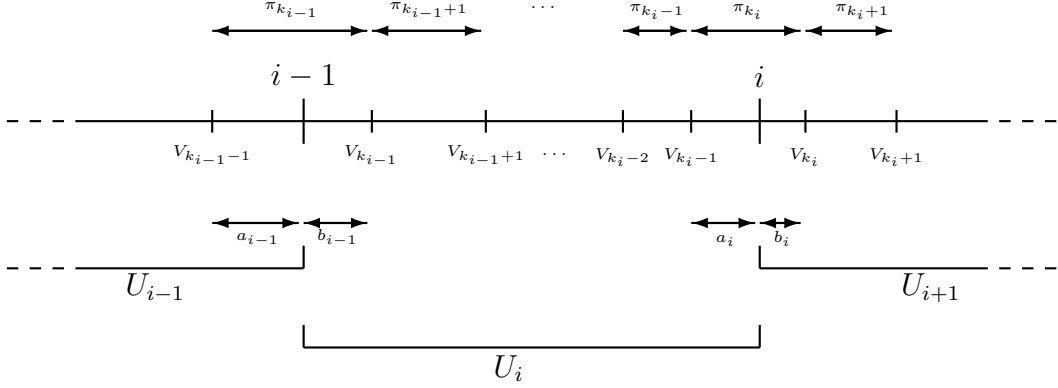

\noindent The microstratum $U_i,~i=1,\ldots,n,$ is defined as
    \begin{equation} \label{Ui}
      U_i=\{k \in U;~k_{i-1} \leq k \leq k_i\}.
    \end{equation}
For any unit $k \in U_i$, we note
    \begin{eqnarray} \label{alpha:ik}
    \alpha_{ik} & = & \left\{\begin{array}{ll}
                               b_{i-1} & \textrm{if } k=k_{i-1}, \\
                               \pi_k & \textrm{if } k_{i-1}<k<k_i, \\
                               a_i & \textrm{if } k=k_i,
                             \end{array}
                      \right.
    \end{eqnarray}
and we note $\alpha_i = (b_{i-1},\pi_{k_{i-1}+1},\ldots,\pi_{k_{i}-1},a_i)^{\top}$. We have in particular
    $$\sum_{k \in U_i} \alpha_{ik}=1.$$
To fix ideas, useful quantities for population $U$ are presented in Figure \ref{nota}. The microstrata are overlapping, since one cross-border unit belongs to two adjacent microstrata: the cross-border unit $k_i$ belongs both to the microstratum $U_i$ (with an associated probability $a_i$) and to the microstratum $U_{i+1}$ (with an associated probability $b_i$).

\section{Ordered pivotal sampling} \label{sec:ops}

\noindent Ordered pivotal sampling (Deville and Till\'e, 1998) may be recursively defined as follows. We initialize with $L_0=k_0$. The unit $L_{i-1}$ jumps from microstratum $U_{i-1}$ to microstratum $U_{i}$, with the residual probability $b_{i-1}$. One unit, denoted as $S_i$, is selected among $\{L_{i-1},k_{i-1}+1,\ldots,k_{i}-1\}$ with probabilities proportional to $(b_{i-1},\pi_{k_{i-1}+1},\ldots,\pi_{k_{i}-1})$. The unit $S_i$ then faces $k_i$. One of these two units, denoted as $F_i$, is selected while the other one, denoted as $L_i$, jumps to microstratum $U_{i+1}$ with the residual probability $b_i$. More precisely, we have
    \begin{eqnarray} \label{Fi:Li}
      (F_i,L_i) & = & \left\{\begin{array}{ll}
                               (S_i,k_i) & \textrm{with probability } \frac{1-a_i-b_i}{1-b_i}, \\
                               (k_i,S_i) & \textrm{with probability } \frac{a_i}{1-b_i}.
                             \end{array}
      \right.
    \end{eqnarray}
It can be shown that this sampling algorithm enables to match exactly the set $\pi_U$ of prescribed inclusion probabilities (Deville and Till\'e, 1998).

\begin{prop} \label{prop1}
  We can write
    \begin{eqnarray} \label{eq1:prop1}
      \hat{t}_{y\pi}-t_y = \sum_{i=1}^n \xi_i & \textrm{ where } & \xi_i = \check{y}_{F_i}+b_i \check{y}_{L_i} - \left\{\sum_{k \in U'_{i}} \alpha_{ik} \check{y}_{k}+ b_i \check{y}_{k_i}\right\},
    \end{eqnarray}
  and where $U'_i = \{L_{i-1},k_{i-1}+1,\ldots,k_{i}-1,k_i\}$. Also, we have
    \begin{eqnarray} \label{eq2:prop1}
      E_{\{y_U,L_{i-1}\}}(\xi_i) & = & 0,
    \end{eqnarray}
  and the $\xi_i$'s are not correlated.
\end{prop}

\begin{propproof}
  Using the identity
    \begin{eqnarray*}
      \sum_{k \in U'_{i}} \alpha_{ik} \check{y}_{k} & = & b_{i-1} \check{y}_{L_{i-1}} + \sum_{k_{i-1}<k<k_i} \alpha_{ik} \check{y}_{k} +  a_i \check{y}_{k_i},
    \end{eqnarray*}
  we can write
    \begin{eqnarray} \label{eq1:pprop1}
      \sum_{i=1}^n \xi_i = \sum_{i=1}^n \left(\check{y}_{F_i}+b_i \check{y}_{L_i} - b_{i-1} \check{y}_{L_{i-1}}\right)
                              +\sum_{i=1}^n \left\{\sum_{k_{i-1}<k<k_i} \alpha_{ik} \check{y}_{k}+(a_i+b_i) \check{y}_{k_i}\right\}.
    \end{eqnarray}
  We first consider the first term in the right-hand side of (\ref{eq1:pprop1}), which we can write as
    \begin{eqnarray}
      \sum_{i=1}^n \left(\check{y}_{F_i}+b_i \check{y}_{L_i} - b_{i-1} \check{y}_{L_{i-1}}\right)
        & = & \sum_{i=1}^n \check{y}_{F_i}+ \sum_{i=1}^n (b_i \check{y}_{L_i}-b_{i-1} \check{y}_{L_{i-1}}) \nonumber \\
        & = & \hat{t}_{y\pi}+(b_n \check{y}_{L_n}-b_{0} \check{y}_{L_{0}}) \nonumber \\
        & = & \hat{t}_{y\pi}. \label{eq2:pprop1}
    \end{eqnarray}
  We now consider the second term in the right-hand side of (\ref{eq1:pprop1}), which we can write as
      \begin{eqnarray}
      \sum_{i=1}^n \left\{\sum_{k_{i-1}<k<k_i} \alpha_{ik} \check{y}_{k}+(a_i+b_i) \check{y}_{k_i}\right\}
        & = & \sum_{i=1}^n \sum_{k_{i-1} < k \leq k_i} y_{k} \nonumber \\
        & = & t_y. \label{eq3:pprop1}
    \end{eqnarray}
  By plugging (\ref{eq2:pprop1}) and (\ref{eq3:pprop1}) into (\ref{eq1:pprop1}), we obtain (\ref{eq2:prop1}). Also, from equation (\ref{Fi:Li}), we obtain
    \begin{eqnarray} \label{eq4:pprop1}
      E_{\{y_U,L_{i-1},S_i\}}(\check{y}_{F_i}+b_i \check{y}_{L_i}) & = & (1-a_i) \check{y}_{S_i}+(a_i+b_i) \check{y}_{k_i}
    \end{eqnarray}
  and from the definition of $\xi_i$ and $S_i$, we obtain
    \begin{eqnarray*}
      E_{\{y_U,L_{i-1}\}}(\xi_{i}) & = & E_{\{y_U,L_{i-1}\}} E_{\{y_U,L_{i-1},S_i\}} (\xi_{i}) \\
                               & = & E_{\{y_U,L_{i-1}\}} \left\{(1-a_i) \check{y}_{S_i} - \sum_{k \in U'_{i} \setminus \{k_i\}} \alpha_{ik} \check{y}_{k} \right\}\\
                               & = & 0.
    \end{eqnarray*}
  We now prove that the $\xi_i$'s are not correlated. For $i<j$, we have
    \begin{eqnarray}
      {Cov}_{\{y_U\}}(\xi_i,\xi_j) & = & E_{\{y_U\}} {Cov}_{\{y_U,\xi_i,L_{j-1}\}}(\xi_i,\xi_j) \nonumber \\
                                   & + & {Cov}_{\{y_U\}}[{E}_{\{y_U,\xi_i,L_{j-1}\}}(\xi_i),{E}_{\{y_U,\xi_i,L_{j-1}\}}(\xi_j)] \nonumber \\
                       & = & {Cov}_{\{y_U\}}[\xi_i,{E}_{\{y_U,\xi_i,L_{j-1}\}}(\xi_j)]. \label{eq5:pprop1}
    \end{eqnarray}
  Also, $\xi_j$ does not depend on $\xi_i$ conditionally on $L_{j-1}$, which leads to
    \begin{eqnarray*}
      {E}_{\{y_U,\xi_i,L_{j-1}\}}(\xi_j) = {E}_{\{y_U,L_{j-1}\}}(\xi_j)=0
    \end{eqnarray*}
  from equation (\ref{eq2:prop1}). This completes the proof.
\end{propproof}

\begin{prop} \label{prop2}
  We have
    \begin{eqnarray} \label{eq1:prop2}
      \sum_{i=1}^n E_{\{y_U\}}(\xi_{i}^4) & \leq & 16 \left\{2+\frac{1}{1-\pi_M}\right\} \left\{\sum_{l \in U} \pi_l \left(\check{y}_{l}-\frac{t_y}{n}\right)^4 \right\}.
    \end{eqnarray}
\end{prop}

\begin{propproof}
We can write
    \begin{eqnarray} \label{eq1:pprop2}
      \sum_{i=1}^n E_{\{y_U\}}(\xi_{i}^4) & = & E_{\{y_U\}}\left[\sum_{i=1}^n E_{\{y_U,L_{i-1}\}}(\xi_{i}^4)\right].
    \end{eqnarray}
Applying equation (\ref{e2:lem2}) in Lemma \ref{lem2} (see Appendix \ref{append:prop2}) with $\phi(x)=x^4$, we obtain
    \begin{eqnarray} \label{eq2:pprop2}
      \sum_{i=1}^n E_{\{y_U\}}(\xi_{i}^4) & \leq & E_{\{y_U\}}\left[\sum_{i=1}^n \sum_{k,l \in U'_i} \alpha_{ik} \alpha_{il} (\check{y}_{l}-\check{y}_{k})^4\right] \nonumber \\
      & = & E_{\{y_U\}}\left[\sum_{i=1}^n \sum_{k,l \in U'_i} \alpha_{ik} \alpha_{il} \left\{\left(\check{y}_{l}-\frac{t_y}{n}\right)-\left(\check{y}_{k}-\frac{t_y}{n}\right)\right\}^4\right] \nonumber \\
      & \leq & 16 E_{\{y_U\}}\left[\sum_{i=1}^n \sum_{k,l \in U'_i} \alpha_{ik} \alpha_{il} \left(\check{y}_{k}-\frac{t_y}{n}\right)^4\right] \nonumber \\
      & = & 16 E_{\{y_U\}}\left[\sum_{i=1}^n \sum_{k \in U'_i} \alpha_{ik} \left(\check{y}_{k}-\frac{t_y}{n}\right)^4\right].
    \end{eqnarray}
We note $z_k=\left(\check{y}_{k}-\frac{t_y}{n}\right)^4$. From equation (\ref{eq2:pprop2}), we obtain
    \begin{eqnarray}
      \sum_{i=1}^n E_{\{y_U\}}(\xi_{i}^4) & \leq & 16 E_{\{y_U\}}\left[\sum_{i=1}^n \left\{b_{i-1} z_{L_{i-1}} + \sum_{k \in U_i} \alpha_{ik} z_k \right\}\right] \nonumber \\
                                    & \leq & 16 \left\{ E_{\{y_U\}}\left(\sum_{i=1}^n b_{i-1} z_{L_{i-1}} \right) + \sum_{l \in U} \pi_{l} z_l \right\}. \label{eq3:pprop2}
    \end{eqnarray}
By applying Lemma \ref{lem3} in Appendix \ref{append:prop2}, equation (\ref{eq3:pprop2}) leads to (\ref{eq1:prop2}).
\end{propproof}

\begin{prop} \label{prop3}
  We have
    \begin{eqnarray} \label{eq1:prop3}
      V_{\{y_U\}}(\hat{t}_{y\pi}) & \geq & \left\{1-\pi_M\right\}^2 \left\{\sum_{i=1}^n \sum_{k \in U_i} \alpha_{ik} \left(\check{y}_{k}-\sum_{l \in U_i} \alpha_{il} \check{y}_{l}\right)^2 \right\}.
    \end{eqnarray}
\end{prop}

\begin{propproof}
From Proposition \ref{prop1}, the $\xi_i$'s are not correlated and $E_{\{y_U,L_{i-1}\}}(\xi_i)=0$, which leads to
    \begin{eqnarray}
      V_{\{y_U\}}(\hat{t}_{y\pi}) = \sum_{i=1}^n V_{\{y_U\}}(\xi_i) = \sum_{i=1}^n E_{\{y_U\}} V_{\{y_U,L_{i-1}\}}(\xi_i). \label{eq1:pprop3}
    \end{eqnarray}
From Lemma \ref{lem4} in Appendix \ref{append:prop3}, we obtain
    \begin{eqnarray}
      V_{\{y_U,L_{i-1}\}}(\xi_i) & \geq & (1-a_i-b_i) \left\{\sum_{k<l \in U'_i \setminus \{k_i\}} \alpha_{ik} \alpha_{il} (\check{y}_{k}-\check{y}_{l})^2 \right. \nonumber \\
      & + & \left. \sum_{k \in U'_i \setminus \{k_i\}} \alpha_{ik} a_i (\check{y}_{k}-\check{y}_{k_i})^2\right\} \nonumber \\
      & = & (1-a_i-b_i) \sum_{k,l \in U'_i} \alpha_{ik} \alpha_{il} (\check{y}_{k}-\check{y}_{l})^2 \nonumber \\
      & = & (1-a_i-b_i) \sum_{k \in U'_i} \alpha_{ik} \left(\check{y}_{k}-\sum_{l \in U'_i} \alpha_{il} \check{y}_{l}\right)^2. \label{eq2:pprop3}
    \end{eqnarray}
Since $L_{i-1}=k_{i-1}$ with probability $\frac{1-a_{i-1}-b_{i-1}}{1-b_{i-1}}$, we obtain from (\ref{eq2:pprop3})
    \begin{eqnarray}
      E_{\{y_U\}} V_{\{y_U,L_{i-1}\}}(\xi_i) & \geq & \frac{(1-a_{i-1}-b_{i-1})(1-a_i-b_i)}{1-b_{i-1}} \sum_{k \in U_i} \alpha_{ik} \left(\check{y}_{k}-\sum_{l \in U_i} \alpha_{il} \check{y}_{l}\right)^2 \nonumber \\
                                             & \geq & \left\{1-\max_{k \in U} \pi_k\right\}^2 \left\{\sum_{k \in U_i} \alpha_{ik} (\check{y}_{k}-\sum_{l \in U_i} \alpha_{il} \check{y}_{l})^2\right\}. \label{eq3:pprop3}
    \end{eqnarray}
From (\ref{eq1:pprop3}) and (\ref{eq3:pprop3}), we obtain (\ref{eq1:prop3}).
\end{propproof}

\begin{prop} \label{prop4}
  We have
    \begin{eqnarray} \label{eq1:prop4}
      V_{\{y_U\}}\left\{\sum_{i=1}^n V_{\{y_U,L_{i-1}\}}(\xi_i)\right\}
      \leq 8\left\{3+\frac{2}{1-\pi_M}\right\}\left\{2+\frac{1}{1-\pi_M}\right\} \sum_{k \in U} \pi_k\left(\check{y}_{k}-\frac{t_y}{n}\right)^4.
    \end{eqnarray}
\end{prop}

The proof of Proposition \ref{prop4} is lengthy, and is thus given in Appendix \ref{append:prop4}.

\section{A design-based central-limit theorem} \label{sec:tcl}

In this Section, we follow the approach in Ohlsson~(1986) and prove the asymptotic normality of the HT-estimator under ordered pivotal sampling by applying the martingale central-limit theorem. A version of this Theorem is reminded in Appendix \ref{appen:mclt}. Let $\sigma(X)$ denote the $\sigma$-field generated by some random variable $X$. We introduce the following $\sigma$-fields:
    \begin{eqnarray} \label{filt:1}
      \mathcal{F}_{0} & = & \sigma(y_{U}), \nonumber \\
      \mathcal{F}_{i} & = & \sigma(y_{U},S_{1},F_{1},L_{1},\ldots,S_{i},F_{i},L_{i}) \textrm{ for } i=1,\ldots,n.
    \end{eqnarray}
We note
    \begin{eqnarray} \label{martingale:1}
      \eta_{i} & = & \frac{\xi_{i}}{\sqrt{V_{\{y_U\}}(\hat{t}_{y\pi})}}
    \end{eqnarray}
where $\xi_{i}$ is defined in (\ref{eq1:prop1}). It follows from Proposition \ref{prop1} that $\{\eta_{i};~i=1,\ldots,n\}$ is a martingale difference sequence with respect to the filtration $\{\mathcal{F}_{i};~i=0,\ldots,n\}$. We make the following assumptions:
\begin{itemize}
  \item[H1:] There exists some constant $f_1<1$ such that for any $k \in U$:
    \begin{eqnarray} \label{eq1:H1}
      \pi_{k} & \leq & f_1.
    \end{eqnarray}
  \item[H2:] There exists some constant $C_1$ such that:
    \begin{eqnarray} \label{eq1:H2}
      \sum_{k \in U} \pi_{k} \left(\check{y}_{k}-\frac{t_{y}}{n}\right)^4 & \leq & C_1 \frac{N^4}{n^3}.
    \end{eqnarray}
  \item[H3:] There exists some constant $C_2>0$ such that:
    \begin{eqnarray} \label{eq1:H3}
      \sum_{i=1}^{n} \sum_{k \in U_i} \alpha_{ik} \left(\check{y}_{k}-\sum_{l \in U_i} \alpha_{il} \check{y}_{l}\right)^2 & \geq & C_2 \frac{N^2}{n}.
    \end{eqnarray}
\end{itemize}

It assumed in (H1) that the first-order inclusion probabilities are bounded away from $1$. This is not a severe restriction in practice, since some unit with an inclusion probability equal to $1$ is automatically surveyed, and is thus not involved in the selection process. The condition (H2) will hold in particular if the variable $y$ has a finite centered moment or order $4$, and if there exists some constants $D_1,D_2>0$ such that for any $k \in U,~D_1 \leq N n^{-1} \pi_k \leq D_2$. Assumption (H3) requires that the dispersion within the micro-strata does not vanish.

\begin{theo} \label{theo1}
  Suppose that the sample $S$ is selected by means of ordered pivotal sampling, and that assumptions (H1)-(H3) hold. Then
    \begin{eqnarray}
      \frac{\hat{t}_{y\pi}-t_y}{\sqrt{V_{\{y_U\}}(\hat{t}_{y\pi})}} & \underset{\mathcal{L}_p}{\longrightarrow} & \mathcal{N}(0,1),
    \end{eqnarray}
  where $\underset{\mathcal{L}_p}{\longrightarrow}$ stands for the convergence in distribution with respect to the randomization associated to the sampling design $p(\cdot)$.
\end{theo}

\begin{theoproof}
  It is sufficient to prove that the sequence $\{\eta_{i};~i=1,\ldots,n\}$ fulfills conditions (a) and (b) in Proposition \ref{appen:theo:mclt} (see Appendix \ref{appen:mclt}). Condition (a) will follow with $\delta=2$ from
    \begin{eqnarray} \label{ptheo1:eq1}
      \sum_{i=1}^{n} E_{\{y_U\}}(\eta_{i}^4) & = & O(n^{-1}).
    \end{eqnarray}
  Also, since $E_{\{y_U\}}\left\{\sum_{i=1}^n V_{\mathcal{F}_{i}}(\eta_{i})\right\} = 1$, condition (b) will follow from
    \begin{eqnarray} \label{ptheo1:eq2}
      V_{\{y_U\}}\left\{\sum_{i=1}^n V_{\mathcal{F}_{i}}(\eta_{i})\right\} & = & O(n^{-1}).
    \end{eqnarray}
  From Propositions \ref{prop2} and \ref{prop3}, we obtain
    \begin{eqnarray*}
      \sum_{i=1}^{n} E_{\{y_U\}}(\eta_{i}^4) & \leq & 16 \frac{\left\{2+\frac{1}{1-\pi_M}\right\} \left\{\sum_{l \in U} \pi_l \left(\check{y}_{l}-\frac{t_y}{n}\right)^4 \right\}}{\left\{1-\pi_M\right\}^4 \left\{\sum_{i=1}^n \sum_{k \in U_i} \alpha_{ik} \left(\check{y}_{k}-\sum_{l \in U_i} \alpha_{il} \check{y}_{l}\right)^2 \right\}^2}
    \end{eqnarray*}
  so that from Assumptions $(H1)-(H3)$, we obtain (\ref{ptheo1:eq1}). Also, from Propositions \ref{prop3} and \ref{prop4}, we obtain
    \begin{eqnarray} \label{ptheo1:eq3}
     V_{\{y_U\}}\left\{\sum_{i=1}^n V_{\mathcal{F}_{i}}(\eta_{i})\right\} \leq 8\frac{\left\{3+\frac{2}{1-\pi_M}\right\}\left\{2+\frac{1}{1-\pi_M}\right\} \sum_{k \in U} \pi_k\left(\check{y}_{k}-\frac{t_y}{n}\right)^4}{\left\{1-\pi_M\right\}^4 \left\{\sum_{i=1}^n \sum_{k \in U_i} \alpha_{ik} \left(\check{y}_{k}-\sum_{l \in U_i} \alpha_{il} \check{y}_{l}\right)^2 \right\}^2}.
    \end{eqnarray}
  From Assumptions $(H1)-(H3)$, we obtain (\ref{ptheo1:eq2}) which completes the proof.
\end{theoproof}

\noindent From Theorem \ref{theo1}, an approximate two-sided $100(1-2\alpha) \% $ confidence interval for $t_y$ is thus given by
    \begin{eqnarray} \label{conf:int:p}
      \left[\hat{t}_{y\pi} \pm u_{1-\alpha} \sqrt{V_{\{y_U\}}(\hat{t}_{y\pi})} \right]
    \end{eqnarray}
with $u_{1-\alpha}$ the quantile of order $1-\alpha$ of the standard normal distribution. In practice, the variance of the HT-estimator in (\ref{conf:int:p}) needs to be replaced with a variance estimator. A customary choice for a fixed-size sampling design is the Sen-Yates-Grundy variance estimator
    \begin{eqnarray} \label{syg:var:est}
      \hat{V}_{SYG} & = & \frac{1}{2} \sum_{k \neq l \in S} \frac{\pi_k \pi_l-\pi_{kl}}{\pi_{kl}}(\check{y}_{k}-\check{y}_{l})^2,
    \end{eqnarray}
with $\pi_{kl}$ the probability that units $k$ and $l$ are selected jointly in the sample. These probabilities can be computed exactly for ordered pivotal sampling (Chauvet, 2012), but many of them are usually equal to $0$ which results in a variance estimator biased downwards. A necessary and sufficient condition for all second-order inclusion probabilities to be strictly positive is that all micro-strata contain at most one non cross-border unit, which is unlikely to hold when the first-order inclusion probabilities are small.

\section{A model-assisted central limit theorem} \label{sec:model}

\noindent To handle this problem, we appeal to some super-population model. We assume that the first-order inclusion probabilities $\pi_k$ are deterministic quantities, possibly defined proportionally to some auxiliary variable $x_k$, and that the values of the $y$-variable are generated according to
   \begin{eqnarray} \label{superpop:model}
      m:y_{k} = \beta~\pi_{k}+ \pi_{k} \epsilon_{k} & \textrm{ for $k \in U$, where } & \left\{\begin{array}{l}
                                                                                                      E(\epsilon_{k})=0, \\
                                                                                                      {V}(\epsilon_{k})=\sigma^2,
                                                                                                     \end{array} \right.
   \end{eqnarray}
and where $\beta$ and $\sigma$ are unknown parameters. The $\epsilon_k$'s need not be independent and identically distributed, and in applications like in spatial sampling they are typically not. We note $\epsilon_U=(\epsilon_1,\ldots,\epsilon_N)^{\top}$, and
    \begin{eqnarray} \label{martingale:2}
      \eta'_{i} & = & \frac{\xi_{i}}{\sqrt{V(\hat{t}_{y\pi}-t_y)}}.
    \end{eqnarray}
We have
    \begin{eqnarray} \label{antic:var}
    V(\hat{t}_{y\pi}-t_y) & = & E V_{\{y_U\}}(\hat{t}_{y\pi}-t_y).
    \end{eqnarray}
Also, it follows from Proposition \ref{prop1} that $\{\eta'_{i};~i=1,\ldots,n\}$ is a martingale difference sequence with respect to the filtration $\{\mathcal{F}_{i};~i=0,\ldots,n\}$. We consider the following assumptions:
\begin{itemize}
  \item[H2b:] There exists some constant $C_3$ such that:
    \begin{eqnarray} \label{eq1:H2b}
      \sum_{k \in U} \pi_{k} E(\epsilon_k^4) & \leq & C_3 \frac{N^4}{n^3}.
    \end{eqnarray}
  \item[H3b:] There exists some constant $C_4>0$ such that:
    \begin{eqnarray} \label{eq1:H3b}
     \sum_{i=1}^{n} \sum_{k,l \in U_i} \alpha_{ik} \alpha_{il} E(\epsilon_k -\epsilon_l)^2 & \geq & C_4 \frac{N^2}{n}.
    \end{eqnarray}
  \item[H4:] There exists some constant $C_5$ such that:
    \begin{flalign} \label{eq1:H4}
      & \sum_{k \neq l \in U} \pi_k(1-\pi_k) \pi_l(1-\pi_l) {Cov}^{+}(\epsilon_k^2,\epsilon_l^2) \leq C_5 \frac{N^3}{n^2}, & \nonumber \\
      & \sum_{k \neq l \neq i \in U} (\pi_{k} \pi_{l}-\pi_{kl}) \pi_i(1-\pi_i) {Cov}^{+}(\epsilon_k \epsilon_l,\epsilon_i^2) \leq C_5 \frac{N^3}{n^2}, & \\
      & \sum_{k \neq l \neq i \neq j \in U} (\pi_{k} \pi_{l}-\pi_{kl}) (\pi_{i} \pi_{j}-\pi_{ij}){Cov}(\epsilon_k \epsilon_l,\epsilon_i \epsilon_j) \leq C_5 \frac{N^3}{n^2}, & \nonumber
    \end{flalign}
  where for two random variables $X$ and $Y$, we note ${Cov}^{+}(X,Y)=Max\{0,{Cov}(X,Y)\}$.
\end{itemize}
Assumptions (H2b) and (H3b) are model counterparts for assumptions (H2) and (H3). The assumption (H4) requires that there is no long-range dependence in the population $U$. This assumption will hold automatically if the $\epsilon_k$'s are independent, since in such case the left-hand side in each of the three lines of (\ref{eq1:H4}) is equal to $0$.

\begin{theo} \label{theo2}
  Suppose that the super-population model $m$ in (\ref{superpop:model}) holds. Suppose that the sample $S$ is selected by means of ordered pivotal sampling, and that assumptions (H1), (H2b), (H3b) and (H4) hold. Then
    \begin{eqnarray}
      \frac{\hat{t}_{y\pi}-t_y}{\sqrt{V(\hat{t}_{y\pi}-t_y)}} & \underset{\mathcal{L}_{mp}}{\longrightarrow} & \mathcal{N}(0,1),
    \end{eqnarray}
  where $\underset{\mathcal{L}_{mp}}{\longrightarrow}$ stands for the convergence in distribution with respect to the randomization associated to both the super-population model $m$ and the sampling design $p(\cdot)$.
\end{theo}

\begin{theoproof}
  It is sufficient to prove that the sequence $\{\eta'_{i};~i=1,\ldots,n\}$ fulfills conditions (a) and (b) in Proposition \ref{appen:theo:mclt}. Condition (a) will follow with $\delta=2$ from
    \begin{eqnarray} \label{ptheo2:eq1}
      \sum_{i=1}^{n} E\{(\eta'_{i})^4\} & = & O(n^{-1}).
    \end{eqnarray}
  From Proposition \ref{prop2}, we obtain after some algebra that there exists some constant $C_6$ such that
    \begin{eqnarray} \label{ptheo2:eq2}
      \sum_{i=1}^n E(\xi_{i}^4) & \leq & C_6 \left\{2+\frac{1}{1-\pi_M}\right\} \left\{\sum_{l \in U} \pi_l E(\epsilon_k)^4 \right\}.
    \end{eqnarray}
  From Proposition \ref{prop3}, we obtain
    \begin{eqnarray} \label{ptheo2:eq3}
      E V_{\{y_U\}}(\hat{t}_{y\pi}) & \geq & \left\{1-\pi_M\right\}^2 \sum_{i=1}^n \sum_{k,l \in U_i} \alpha_{ik} \alpha_{il} E(\epsilon_k -\epsilon_l)^2.
    \end{eqnarray}
  From equations (\ref{ptheo2:eq2}) and (\ref{ptheo2:eq3}), and from Assumptions (H1), (H2b) and (H3b), we obtain (\ref{ptheo2:eq1}). \\

  We now turn to condition (b), which will follow from
    \begin{eqnarray} \label{ptheo2:eq4}
      V \left\{\sum_{i=1}^n V_{\mathcal{F}_{i}}(\eta'_{i})\right\} & = & O(n^{-1}).
    \end{eqnarray}
  We have
    \begin{eqnarray} \label{ptheo2:eq5}
      V \left\{\sum_{i=1}^n V_{\mathcal{F}_{i}}(\eta'_{i})\right\} = E\left[V_{\{y_U\}}\left\{\sum_{i=1}^n V_{\mathcal{F}_{i}}(\eta'_{i})\right\}\right] + \frac{V\left[V_{\{y_U\}}(\hat{t}_{y\pi})\right]}{\left\{V(\hat{t}_{y\pi}-t_y)\right\}^2}.
    \end{eqnarray}
  Making use of Proposition \ref{prop4} and of Assumptions (H1), (H2b) and (H3b), it is easily proved that the first term in the right-hand side of (\ref{ptheo2:eq5}) is $O(n^{-1})$. The proof (omitted) is similar to that for equation (\ref{ptheo2:eq1}). Also we have
    \begin{eqnarray} \label{ptheo2:eq6}
      V\left[V_{\{y_U\}}(\hat{t}_{y\pi})\right] & = & V\left[\frac{1}{2} \sum_{k \neq l \in U} (\pi_{k} \pi_{l}-\pi_{kl})(\check{y}_k-\check{y}_l)^2 \right] \\
                                                & = & \frac{1}{4} \sum_{k \neq l \in U} \sum_{i \neq j \in U} (\pi_{k} \pi_{l}-\pi_{kl})(\pi_{i} \pi_{j}-\pi_{ij}){Cov}\left\{(\epsilon_k-\epsilon_l)^2,(\epsilon_i-\epsilon_j)^2 \right\}. \nonumber
    \end{eqnarray}
  After some algebra, we obtain
    \begin{eqnarray} \label{ptheo2:eq7}
      V\left[V_{\{y_U\}}(\hat{t}_{y\pi})\right] & \leq & \sum_{k \neq l \in U} \pi_k(1-\pi_k) \pi_l(1-\pi_l) {Cov}^{+}(\epsilon_k^2,\epsilon_l^2) \nonumber \\
      & + & 2 \sum_{k \neq l \neq i \in U} (\pi_{k} \pi_{l}-\pi_{kl}) \pi_i(1-\pi_i) {Cov}^{+}(\epsilon_k \epsilon_l,\epsilon_i^2) \\
      & + & \sum_{k \neq l \neq i \neq j \in U} (\pi_{k} \pi_{l}-\pi_{kl}) (\pi_{i} \pi_{j}-\pi_{ij}){Cov}(\epsilon_k \epsilon_l,\epsilon_i \epsilon_j), \nonumber
    \end{eqnarray}
  which is $O(N^3 n^{-2})$ from Assumption (H4). From equation (\ref{ptheo2:eq3}) and Assumption (H3b), we obtain that the second term in the right-hand side of (\ref{ptheo2:eq5}) is $O(N^{-1})$. This completes the proof.
\end{theoproof}

\noindent From Theorem \ref{theo2}, an approximate two-sided $100(1-2\alpha) \% $ confidence interval for $t_y$ is thus given by
    \begin{eqnarray} \label{conf:int:mp}
      \left[\hat{t}_{y\pi} \pm u_{1-\alpha} \sqrt{V(\hat{t}_{y\pi}-t_y)} \right]
    \end{eqnarray}
with $u_{1-\alpha}$ the quantile of order $1-\alpha$ of the standard normal distribution. To obtain a variance estimator, we first note that under the super-population model $m$:
    \begin{eqnarray} \label{antic:var:2}
      V(\hat{t}_{y\pi}-t_y) & = & \sigma^2 \sum_{k \in U} \pi_k(1-\pi_k) + \sum_{k \neq l \in U} (\pi_{kl}-\pi_k \pi_l) {Cov}(\epsilon_k,\epsilon_l).
    \end{eqnarray}
If the $\epsilon_k$'s are assumed to be independent, the second term in (\ref{antic:var:2}) vanishes which leads to the variance estimator
    \begin{eqnarray} \label{est:var:mp}
      \hat{V}_m = \hat{\sigma}^2 \sum_{k \in U} \pi_k(1-\pi_k) & \textrm{ where } & \hat{\sigma}^2=\frac{1}{2n(n-1)} \sum_{k \in S} (\check{y}_l-\check{y}_k)^2. \end{eqnarray}
This variance estimator is then unbiased for $V(\hat{t}_{y\pi}-t_y)$, provided that $S$ is independent on $\epsilon_U$ (non-informative sampling). Alternatively, a weighted version of $\hat{\sigma}^2$ may be used. If the $\epsilon_k$'s may not be assumed to be independent, like in case of spatial sampling for example, the covariance structure of the $\epsilon_k$'s may be parametrically modeled. Several methods for the estimation of these parameters are possible, see for example Sherman~(2011).

\section{Conclusion} \label{sec:conc}

\noindent In this work, we proved a central-limit theorem for the HT-estimator both under a design-based approach and under a model-assisted approach. The question of variance estimation was only briefly discussed. It would be of both theoretical and practical interest to exhibit variance estimators which are weakly consistent, under reasonable model assumptions, so as to build confidence intervals in practice with approximate appropriate coverage. This goes beyond the goal of the present paper, but is a matter for further research. \\

\noindent An auxiliary variable is needed to define the first-order inclusion probabilities. In practice, other auxiliary variables may be available for all the units in the population, which would enable to consider a more elaborate super-population model. In particular, balanced sampling by means of the cube method could be used to achieve a spatial balance, while controlling the balancing on these auxiliary variables (Grafstr\"om and Till\'e, 2013). Establishing a central-limit theorem for the cube method and an arbitrary set of auxiliary variables is also a matter for further research.

\section*{References}

Br\"and\'en, P., and Jonasson, J. (2011). Negative dependence in sampling. \textit{Scand. J. Stat.} \textbf{39} 830--838. \\

Chauvet, G. (2012). On a characterization of ordered pivotal sampling. \textit{Bernoulli.} \textbf{18} 1320--1340. \\

Deville, J-C., and Till\'e, Y. (1998). Unequal probability sampling without replacement through a splitting method. \textit{Biometrika.} \textbf{85} 89--101. \\

Grafstr\"om, A., and Lundstr\"om, N. L., and Schelin, L. (2012). Spatially balanced sampling through the pivotal method. \textit{Biometrics.} \textbf{68} 514--520. \\

Grafstr\"om, A., and Till\'e, Y. (2013). Doubly balanced spatial sampling with spreading and restitution of auxiliary totals. \textit{Environmetrics.} \textbf{24} 120--131. \\

Helland, I.S. (1982). Central limit theorems for martingales with discrete or continuous time. \textit{Scand. J. Stat.} \textbf{9} 79--94. \\

Isaki, C.T., and Fuller, W.A. (1982). Survey design under the regression super-population model. \textit{J. Amer. Statist. Assoc.} \textbf{77} 89--96. \\

Jonasson, J. (2012). The BK inequality for pivotal sampling a.k.a. the Srinivasan sampling process. \textit{Electron. Commun. Probab.} \textbf{35} 1--6. \\

Nedyalkova, D., and Qualit\'e, L., and Till\'e, Y. (2009). General framework for the rotation of units in repeated survey sampling. \textit{Stat. Neerl.} \textbf{63} 269--293. \\

Ohlsson, E. (1986). Asymptotic normality of the Rao-Hartley-Cochran estimator: an Application of the Martingale CLT. \textit{Scand. J. Stat.} \textbf{13} 17--28. \\

Sherman, M. (2011). \textit{Spatial statistics and spatio-temporal data: covariance functions and directional properties}, Wiley, New-York. \\

Srinivasan, A. (2001). Distributions on level sets with applications to approximation algorithms. \textit{Proceedings of the 42nd IEEE Symposium on the Foundations of Computer Science}. 

Till\'e, Y. (2006). \textit{Sampling algorithms}, Springer, New York.

\appendix

\section{Additional results for Proposition \ref{prop2}} \label{append:prop2}

\begin{lemma} \label{lem1}
  We can alternatively write
    \begin{eqnarray}
      \xi_{i} & = & b_i (\check{y}_{S_i}-\check{y}_{F_i})+\sum_{k \in U'_i} \alpha_{ik} (\check{y}_{F_i}-\check{y}_{k}) \label{eq1:lem1} \\
              & = & (1-b_i) \sum_{k \in U'_i} \alpha_{ik}(\check{y}_{F_i}-\check{y}_{k})+b_i \sum_{k \in U'_i} \alpha_{ik}(\check{y}_{S_i}-\check{y}_{k}). \label{eq2:lem1}
    \end{eqnarray}
\end{lemma}

\begin{lemproof}
  By noting that $\check{y}_{F_i}+\check{y}_{L_i}=\check{y}_{S_i}+\check{y}_{k_i}$, we obtain from the definition of $\xi_{i}$:
    \begin{eqnarray*}
      \xi_{i} & = & \check{y}_{F_i}+b_i (\check{y}_{L_i}-\check{y}_{k_i}) - \sum_{k \in U'_i} \alpha_{ik} \check{y}_{k} \\
              & = & \check{y}_{F_i}+b_i (\check{y}_{S_i}-\check{y}_{F_i}) - \sum_{k \in U'_i} \alpha_{ik} \check{y}_{k} \\
              & = & \sum_{k \in U'_i} \alpha_{ik} \check{y}_{F_i}+b_i (\check{y}_{S_i}-\check{y}_{F_i}) - \sum_{k \in U'_i} \alpha_{ik} \check{y}_{k} \\
              & = & b_i (\check{y}_{S_i}-\check{y}_{F_i})+\sum_{k \in U'_i} \alpha_{ik} (\check{y}_{F_i}-\check{y}_{k}).
    \end{eqnarray*}
  Equation (\ref{eq2:lem1}) follows from (\ref{eq1:lem1}).
\end{lemproof}

\begin{lemma} \label{lem2}
  For any convex function $\phi$, we have:
    \begin{eqnarray}
      \phi(\xi_{i} ) & \leq & (1-b_i) \sum_{k \in U'_i} \alpha_{ik} \phi(\check{y}_{F_i}-\check{y}_{k})+b_i \sum_{k \in U'_i} \alpha_{ik} \phi(\check{y}_{S_i}-\check{y}_{k}), \label{eq1:lem2} \\
      E_{\{L_{i-1}\}} \phi(\xi_{\nu i} ) & \leq & \sum_{k,l \in U'_i} \alpha_{ik} \alpha_{il} \phi(\check{y}_{l}-\check{y}_{k}). \label{e2:lem2}
    \end{eqnarray}
\end{lemma}

\begin{lemproof}
  Equation (\ref{eq1:lem2}) readily follows from equation (\ref{eq2:lem1}) and from the Jensen inequality. For any $k \in U'_i$, we have:
    \begin{eqnarray*}
      E_{\{L_{i-1},S_i\}} \phi(\check{y}_{F_i}-\check{y}_{k}) & = & \frac{1-a_i-b_i}{1-b_i} \phi(\check{y}_{S_i}-\check{y}_{k}) + \frac{a_i}{1-b_i} \phi(\check{y}_{k_i}-\check{y}_{k}),
    \end{eqnarray*}
  and from equation (\ref{eq1:lem2}):
    \begin{eqnarray*}
      E_{\{L_{i-1},S_i\}} \phi(\xi_{\nu i} ) & \leq & (1-a_i) \sum_{k \in U'_i} \alpha_{ik} \phi(\check{y}_{S_i}-\check{y}_{k})+ a_i \sum_{k \in U'_i} \alpha_{ik} \phi(\check{y}_{k_i}-\check{y}_{k}), \\
      E_{\{L_{i-1}\}} \phi(\xi_{\nu i} ) & \leq & \sum_{k \in U'_i} \sum_{l \in U'_i \setminus \{k_i\}} \alpha_{ik} \alpha_{il} \phi(\check{y}_{l}-\check{y}_{k})+ \sum_{k \in U'_i} a_i \alpha_{ik} \phi(\check{y}_{k_i}-\check{y}_{k}), \\
      & = & \sum_{k,l \in U'_i} \alpha_{ik} \alpha_{il} \phi(\check{y}_{l}-\check{y}_{k}).
    \end{eqnarray*}
\end{lemproof}

\begin{lemma} \label{lem2b}
  For any $i \leq j=1,\ldots,n$, we note $c_i=\frac{a_i b_i}{(1-a_i)(1-b_i)}$, and
    \begin{eqnarray} \label{eq1:lem2b}
    c(i,j) & = & \left\{\begin{array}{ll}
                        1 & \textrm{if } j=i, \\
                        \prod_{l=i}^{j-1} c_l & \textrm{if } j>i.
                        \end{array} \right.
    \end{eqnarray}
  We have
    \begin{eqnarray} \label{eq2:lem2b}
      \sum_{j=i}^n c(i,j) \leq 1+\frac{1}{1-\pi_M} & \textrm{ and } & \sum_{i=1}^j c(i,j) \leq 1+\frac{1}{1-\pi_M}.
    \end{eqnarray}
\end{lemma}

\begin{lemproof}
For any $l=1,\ldots,n$, we have
    \begin{eqnarray}
      c_l = \frac{a_l b_l}{(1-\pi_{k_l})+a_l b_l} \leq \frac{a_l b_l}{(1-\pi_M)+a_l b_l} \leq \frac{1}{2-\pi_M} \equiv c. \label{eq1:plem3b}
    \end{eqnarray}
For any $j \geq i$, this leads to $c(i,j) \leq c^{j-i}$. We obtain
    \begin{eqnarray*}
      \sum_{j=i}^n c(i,j) \leq \sum_{j=i}^n c^{j-i} \leq \frac{1}{1-c} = 1+\frac{1}{1-\pi_M}
    \end{eqnarray*}
and
    \begin{eqnarray*}
      \sum_{i=1}^j c(i,j) \leq \sum_{i=1}^j c^{j-i} \leq \frac{1}{1-c} = 1+\frac{1}{1-\pi_M}.
    \end{eqnarray*}
\end{lemproof}

\begin{lemma} \label{lem3}
 We note $z$ for some variable of interest, and $z_l^+=max(z_l,0)$. Then
    \begin{eqnarray} \label{eq1:lem3}
      E\left(\sum_{j=1}^{n-1} b_j z_{L_j}\right) & \leq & \left\{1+\frac{1}{1-\pi_M}\right\} \sum_{l \in U} \pi_l z_l^+.
    \end{eqnarray}
\end{lemma}

\begin{lemproof}
  We can write
    \begin{eqnarray} \label{eq1:plem3}
      \sum_{j=1}^{n-1} b_j z_{L_j} = \sum_{l \in U} W_l z_l & \textrm{ where } & W_l = \sum_{j=1}^{n-1} b_j 1(L_j=l).
    \end{eqnarray}
  Take any unit $l \in U$, and let $i$ be the smallest integer for which we may have $L_i=l$ so that
    \begin{eqnarray} \label{eq2:plem3}
      W_l = \sum_{j=i}^{n-1} b_j 1(L_j=l).
    \end{eqnarray}
  For any $j>i$, since $\{L_j=l\}$ implies $\{L_{j-1}=l\}$, we have
    \begin{eqnarray}
      b_j Pr(L_j=l) & = & b_j Pr(L_j=l|L_{j-1}=l) Pr(L_{j-1}=l) \nonumber \\
                    & = & b_j \left(\frac{b_{j-1}}{1-a_j} \frac{a_{j}}{1-b_j} \right) Pr(L_{j-1}=l) \nonumber \\
                    & = & c_j \{b_{j-1} Pr(L_{j-1}=l)\}. \label{eq3:plem3}
    \end{eqnarray}
  Making use of (\ref{eq2:plem3}), we obtain
    \begin{eqnarray}
      E(W_l) & = & \sum_{j=i}^{n-1} b_j Pr(L_j=l) \nonumber \\
             & = & \{b_i Pr(L_i=l)\} \sum_{j=i}^{n} c(i,j) \nonumber \\
             & \leq & \left\{1+\frac{1}{1-\pi_M}\right\} \{b_i Pr(L_i=l)\}, \label{eq4:plem3}
    \end{eqnarray}
  where the last inequality follows from equation (\ref{eq2:lem2b}) in Lemma \ref{lem2b}. There are two cases for unit $l$. If $k_{i-1}<l<k_i$ (non cross-border unit), we have
    \begin{eqnarray}
      b_i Pr(L_i=l) = b_i \frac{\pi_l}{1-a_i} \frac{a_i}{1-b_i} = \pi_l \frac{a_i b_i}{(1-a_i)(1-b_i)} \leq \pi_l. \label{eq8:plem3}
    \end{eqnarray}
   If $l=k_{i}$ (cross-border unit), we have
    \begin{eqnarray}
      b_i Pr(L_i=l) = b_i \frac{1-a_i-b_i}{1-b_i} = \pi_l \frac{b_i (1-\pi_l)}{\pi_l(1-b_i)} \leq \pi_l. \label{eq9:plem3}
    \end{eqnarray}
   By plugging (\ref{eq8:plem3}) or (\ref{eq9:plem3}) in (\ref{eq4:plem3}), and using equation (\ref{eq1:plem3}), we obtain (\ref{eq1:lem3}).
\end{lemproof}

\section{Additional results for Proposition \ref{prop3}} \label{append:prop3}

\begin{lemma} \label{lem4}
  We have
    \begin{eqnarray} \label{eq1:lem4}
      V_{\{L_{i-1}\}}(\xi_i) & = & \sum_{k<l \in U'_i \setminus \{k_i\}} \alpha_{ik} \alpha_{il} (\check{y}_{k}-\check{y}_{l})^2 \nonumber \\
                                           & + & \frac{1-a_i-b_i}{1-a_i} \sum_{k \in U'_i \setminus \{k_i\}} \alpha_{ik} a_i (\check{y}_{k}-\check{y}_{k_i})^2.
    \end{eqnarray}
\end{lemma}

\begin{lemproof}
  We can write
    \begin{eqnarray} \label{eq1:plem4}
      V_{\{L_{i-1}\}}(\xi_i) & = & V_{\{L_{i-1}\}} E_{\{L_{i-1},S_i\}}(\xi_i)+E_{\{L_{i-1}\}} V_{\{L_{i-1},S_i\}}(\xi_i).
    \end{eqnarray}
  We consider the first term in the right-hand side of (\ref{eq1:plem4}). From the definition of $F_i$ and $L_i$, we obtain
    \begin{eqnarray}
      E_{\{L_{i-1},S_i\}}(\xi_i) & = & (1-a_i) \check{y}_{S_i}-\sum_{k \in U'_i \setminus \{k_i\}} \alpha_{ik} \check{y}_{k}, \nonumber \\
      V_{\{L_{i-1}\}} E_{\{L_{i-1},S_i\}}(\xi_i) & = & (1-a_i) \sum_{k \in U'_i \setminus \{k_i\}} \alpha_{ik} \left(\check{y}_{k}- \sum_{l \in U'_i \setminus \{k_i\}} \frac{\alpha_{il}}{1-a_i} \check{y}_{l}\right)^2 \nonumber \\
      & = & \sum_{k<l \in U'_i \setminus \{k_i\}} \alpha_{ik} \alpha_{il} (\check{y}_{k}-\check{y}_{l})^2. \label{eq2:plem4}
    \end{eqnarray}
  We now consider the second term in the right-hand side of (\ref{eq1:plem4}). After some algebra, we obtain
    \begin{eqnarray}
      V_{\{L_{i-1},S_i\}}(\xi_i) & = & a_i(1-a_i-b_i) (\check{y}_{S_i}-\check{y}_{k_i})^2, \nonumber \\
      E_{\{L_{i-1}\}} V_{\{L_{i-1},S_i\}}(\xi_i) & = & \frac{1-a_i-b_i}{1-a_i} \sum_{k \in U'_i \setminus \{k_i\}} \alpha_{ik} a_i (\check{y}_{k}-\check{y}_{k_i})^2. \label{eq3:plem4}
    \end{eqnarray}
  From (\ref{eq2:plem4}) and (\ref{eq3:plem4}), we obtain (\ref{eq1:lem4}).
\end{lemproof}

\newpage
\section{Proof of Proposition \ref{prop4}} \label{append:prop4}

From equation (\ref{eq1:lem4}) in Lemma \ref{lem4}, we can write
    \begin{eqnarray} \label{eq1:pprop4}
    V_{\{L_{i-1}\}}(\xi_i) & = & \sum_{k \in U_i} b_{i-1} \beta_{ik} (\check{y}_{L_{i-1}}-\check{y}_{k})^2 + C_{1i},
    \end{eqnarray}
where
    \begin{eqnarray} \label{eq2:pprop4}
      \beta_{ik} & = & \left\{ \begin{array}{ll}
                                 0 & \textrm{ if } k=k_{i-1}, \\
                                 \alpha_{ik} & \textrm{ if } k_{i-1}<k<k_i, \\
                                 \frac{a_i(1-a_i-b_i)}{1-b_i} & \textrm{ if } k=k_{i}
                               \end{array} \right.
    \end{eqnarray}
and where $C_{1i}$ is non-random. We note
    \begin{eqnarray} \label{eq3:pprop4}
      \bar{y}_{\beta i} & = & \sum_{k \in U_i} \beta'_{ik} \check{y}_{k}
    \end{eqnarray}
where $\beta'_{ik}=\beta_{ik}/\beta_{i}$, and $\beta_i=\sum_{k \in U_i} \beta_{ik}$. We have
    \begin{eqnarray} \label{eq4:pprop4}
      \sum_{k \in U_i} b_{i-1} \beta_{ik} (\check{y}_{L_{i-1}}-\check{y}_{k})^2 & = & \beta_i (\check{y}_{L_{i-1}}-\bar{y}_{\beta i})^2+C_{2i},
    \end{eqnarray}
where $C_{2i}$ is non-random. From equations (\ref{eq1:pprop4}) and (\ref{eq4:pprop4}), we obtain
    \begin{eqnarray} \label{eq5:pprop4}
      V\left\{\sum_{i=1}^n V_{\{L_{i-1}\}}(\xi_i)\right\} & = & V\left\{\sum_{i=1}^n b_{i-1} \beta_i (\check{y}_{L_{i-1}}-\bar{y}_{\beta i})^2 \right\}.
    \end{eqnarray}

From equation (\ref{eq5:pprop4}), we can write $V\left\{\sum_{i=1}^n V_{\{L_{i-1}\}}(\xi_i)\right\}=V_1+V_2$, where
    \begin{eqnarray}
      V_1 & = & \sum_{i=1}^n b_{i-1}^2 \beta_i^2 V[(\check{y}_{L_{i-1}}-\bar{y}_{\beta i})^2], \label{eq6:pprop4}\\
      V_2 & = & 2 \sum_{i<j=1}^n \beta_i \beta_j Cov[b_{i-1}(\check{y}_{L_{i-1}}-\bar{y}_{\beta i})^2,b_{j-1}(\check{y}_{L_{j-1}}-\bar{y}_{\beta j})^2]. \label{eq7:pprop4}
    \end{eqnarray}
We first consider $V_1$, for which we have
    \begin{eqnarray}
      V_1 & \leq & \sum_{i=1}^n b_{i-1}^2 \beta_i^2 E(\check{y}_{L_{i-1}}-\bar{y}_{\beta i})^4 \nonumber \\
          & \leq & 8 \sum_{i=1}^n b_{i-1}^2 \beta_i^2 \left\{E\left(\check{y}_{L_{i-1}}-\frac{t_y}{n}\right)^4+\left(\bar{y}_{\beta i}-\frac{t_y}{n}\right)^4\right\} \nonumber \\
          & \leq & 8 \sum_{i=1}^n b_{i-1} E\left(\check{y}_{L_{i-1}}-\frac{t_y}{n}\right)^4+8 \sum_{i=1}^n \beta_i \left(\bar{y}_{\beta i}-\frac{t_y}{n}\right)^4. \label{eq8:pprop4}
    \end{eqnarray}
We now consider $V_2$, for which we have
    \begin{flalign}
      & Cov[b_{i-1}(\check{y}_{L_{i-1}}-\bar{y}_{\beta i})^2,b_{j-1}(\check{y}_{L_{j-1}}-\bar{y}_{\beta j})^2] & \nonumber \\
      \phantom{sp} & = Cov[b_{i-1}(\check{y}_{L_{i-1}}-\bar{y}_{\beta i})^2,E\left\{b_{j-1}(\check{y}_{L_{j-1}}-\bar{y}_{\beta j})^2|L_{i-1}\right\}] & \nonumber \\
      \phantom{sp} & = Cov[b_{i-1}(\check{y}_{L_{i-1}}-\bar{y}_{\beta i})^2,Pr(L_{j-1}=L_{i-1}) b_{j-1}(\check{y}_{L_{i-1}}-\bar{y}_{\beta j})^2] & \nonumber \\
      \phantom{sp} & = Cov[b_{i-1}(\check{y}_{L_{i-1}}-\bar{y}_{\beta i})^2,b_{i-1} c(i,j) (\check{y}_{L_{i-1}}-\bar{y}_{\beta j})^2] & \nonumber \\
      \phantom{sp} & = b_{i-1}^2 c(i,j) Cov[(\check{y}_{L_{i-1}}-\bar{y}_{\beta i})^2,(\check{y}_{L_{i-1}}-\bar{y}_{\beta j})^2], & \label{eq9:pprop4}
    \end{flalign}
where the fourth line in (\ref{eq9:pprop4}) follows from equation (\ref{eq3:plem3}) in Lemma \ref{lem3}. Also, we have
    \begin{flalign}
      & Cov[(\check{y}_{L_{i-1}}-\bar{y}_{\beta i})^2,(\check{y}_{L_{i-1}}-\bar{y}_{\beta j})^2] & \nonumber \\
      & \leq E[(\check{y}_{L_{i-1}}-\bar{y}_{\beta i})^2(\check{y}_{L_{i-1}}-\bar{y}_{\beta j})^2] & \nonumber \\
      & \leq \frac{1}{2} E[(\check{y}_{L_{i-1}}-\bar{y}_{\beta i})^4]+\frac{1}{2}E[(\check{y}_{L_{i-1}}-\bar{y}_{\beta j})^2] & \nonumber \\
      & \leq 8 E\left(\check{y}_{L_{i-1}}-\frac{t_y}{n}\right)^4 +4 \left(\bar{y}_{\beta i}-\frac{t_y}{n}\right)^4+4 \left(\bar{y}_{\beta j}-\frac{t_y}{n}\right)^4. & \label{eq10:pprop4}
    \end{flalign}
By plugging (\ref{eq9:pprop4}) and (\ref{eq10:pprop4}) in (\ref{eq7:pprop4}), we obtain
    \begin{eqnarray}
      V_2 & \leq & 16 \sum_{i=1}^n b_{i-1} E\left(\check{y}_{L_{i-1}}-\frac{t_y}{n}\right)^4 \sum_{j>i} c(i,j) \nonumber \\
          & + & 8 \sum_{i=1}^n \beta_{i} \left(\bar{y}_{\beta i}-\frac{t_y}{n}\right)^4 \sum_{j>i} c(i,j) \nonumber \\
          & + & 8 \sum_{j=2}^n \beta_{j} \left(\bar{y}_{\beta j}-\frac{t_y}{n}\right)^4 \sum_{i<j} c(i,j). \label{eq11:pprop4}
    \end{eqnarray}
From equation (\ref{eq2:lem2b}) in Lemma \ref{lem2b}, we obtain
    \begin{eqnarray} \label{eq13:pprop4}
      V_2 \leq 16 \left\{1+\frac{1}{1-\pi_M}\right\}\left\{\sum_{i=1}^n b_{i-1} E\left(\check{y}_{L_{i-1}}-\frac{t_y}{n}\right)^4 + \sum_{i=1}^n \beta_{i} \left(\bar{y}_{\beta i}-\frac{t_y}{n}\right)^4 \right\},
    \end{eqnarray}
and from equations (\ref{eq8:pprop4}) and (\ref{eq13:pprop4}), we obtain
    \begin{eqnarray} \label{eq14:pprop4}
      V_1+V_2 \leq 8 \left\{3+\frac{2}{1-\pi_M}\right\}\left\{\sum_{i=1}^n b_{i-1} E\left(\check{y}_{L_{i-1}}-\frac{t_y}{n}\right)^4 + \sum_{i=1}^n \beta_{i} \left(\bar{y}_{\beta i}-\frac{t_y}{n}\right)^4 \right\}.
    \end{eqnarray}
By applying Lemma \ref{lem3}, we have
    \begin{eqnarray} \label{eq15:pprop4}
      \sum_{i=1}^n b_{i-1} E\left(\check{y}_{L_{i-1}}-\frac{t_y}{n}\right)^4 & \leq & \left\{1+\frac{1}{1-\pi_M}\right\} \sum_{k \in U} \pi_k\left(\check{y}_{k}-\frac{t_y}{n}\right)^4.
    \end{eqnarray}
Also, we have
    \begin{eqnarray}
      \sum_{i=1}^n \beta_{i} \left(\bar{y}_{\beta i}-\frac{t_y}{n}\right)^4
      & = & \sum_{i=1}^n \beta_{i} \left\{\sum_{k \in U_i} \beta'_{ik} \left(\check{y}_{k}-\frac{t_y}{n}\right)\right\}^4 \nonumber \\
      & \leq & \sum_{i=1}^n \beta_{i} \sum_{k \in U_i} \beta'_{ik} \left(\check{y}_{k}-\frac{t_y}{n}\right)^4  \nonumber \\
      & = & \sum_{i=1}^n \sum_{k \in U_i} \beta_{ik} \left(\check{y}_{k}-\frac{t_y}{n}\right)^4  \nonumber \\
      & \leq & \sum_{i=1}^n \sum_{k \in U_i} \alpha_{ik} \left(\check{y}_{k}-\frac{t_y}{n}\right)^4  \nonumber \\
      & = & \sum_{k \in U} \pi_k\left(\check{y}_{k}-\frac{t_y}{n}\right)^4, \label{eq16:pprop4}
    \end{eqnarray}
by using the Jensen inequality in the second line. From equations (\ref{eq14:pprop4}), (\ref{eq15:pprop4}) and (\ref{eq16:pprop4}), we obtain (\ref{eq1:prop4}).

\section{A general martingale central-limit theorem} \label{appen:mclt}

\begin{prop}[Helland, 1982; Ohlsson, 1986] \label{appen:theo:mclt}
  For $\nu=1,2,\ldots$, let
    \begin{eqnarray*}
      \{\eta_{\nu i};~i=1,\ldots,n_{\nu}\}
    \end{eqnarray*}
  denote a martingale difference sequence relative to the filtration $\{\mathcal{F}_{\nu i};~1=1,\ldots,n_{\nu}\}$ and let
    \begin{eqnarray}
      S_{\nu} & = & \sum_{i=^1}^{n_{\nu}} \eta_{\nu i}.
    \end{eqnarray}
  Assume that the following conditions are fulfilled:
  \begin{itemize}
    \item[(a)] For some $\delta>0$, we have:
        \begin{eqnarray}
          \sum_{i=1}^{n_{\nu}} E|\eta_{\nu i}|^{2+\delta} & \underset{\nu \to\infty}{\longrightarrow} & 0.
        \end{eqnarray}
    \item[(b)] We have:
        \begin{eqnarray}
          \sum_{i=1}^{n_{\nu}} E_{\{\mathcal{F}_{\nu,i-1}\}} \eta_{\nu i}^{2} & \underset{Pr}{\longrightarrow} & 1.
        \end{eqnarray}
  \end{itemize}
  Then $S_{\nu} \underset{\mathcal{L}}{\longrightarrow} \mathcal{N}(0,1)$, where $\underset{\mathcal{L}}{\longrightarrow}$ stands for the convergence in distribution.
\end{prop}

\end{document}